\documentclass[11pt]{article}
\usepackage[margin=1.1in,a4paper]{geometry}

\usepackage{natbib}
\usepackage{blindtext}
\usepackage{amssymb}    
\usepackage{amsmath}    
\usepackage{comment}
\usepackage{color}  

\usepackage[linecolor=blue!60!,backgroundcolor=blue!10!,textwidth=2.5cm,textsize=scriptsize]{todonotes}
\usepackage[normalem]{ulem}
\input{macro}
\allowdisplaybreaks

\usepackage{array}
\usepackage{booktabs}

\def\D{{\cal D}}

\def\R{{\mathbb R}}

\def\P{{\cal P}}

\def\X{{\cal X}}
\def\Y{{\cal Y}}
\def\PROB {{\mathbb P}}
\def\EXP {{\mathbb E}}
\def\IND{{\mathbb I}}

\def\argmin{\mathop{\rm arg\, min}}

\def\suppY {\mathcal{S}}

\newtheorem{theorem}{Theorem}
\newtheorem{example}{Example}

\usepackage{lastpage}

\title{Metric space valued Fréchet regression}

\author{
L\'aszl\'o Gy\"orfi\textsuperscript{1}\thanks{ Correspondence: \texttt{gyorfi@cs.bme.hu}.} ,\quad
Pierre Humbert\textsuperscript{2}\thanks{Part of this work was carried out while P.H. was working at LPSM of Sorbonne University.} ,\quad
Batiste Le Bars\textsuperscript{3}
}
\date{\small
\textsuperscript{1} Department of Computer Science and Information Theory, Budapest University of Technology and Economics \\
\textsuperscript{2}Université Paris-Saclay, CNRS, Univ Evry, Laboratoire de Mathématiques et Modélisation d'Evry, 91037, Evry-Courcouronnes, France \\
\textsuperscript{3}Univ. Lille, Inria, CNRS, Centrale Lille, UMR 9189, CRIStAL, F-59000 Lille \\
}

\begin{document}

\maketitle

\begin{abstract}
We consider the problem of estimating the Fréchet and conditional Fréchet mean from data taking values in separable metric spaces. Unlike Euclidean spaces, where well-established methods are available, there is no practical estimator that works universally for all metric spaces. Therefore, we introduce a computable estimator for the Fréchet mean based on random quantization techniques and establish its universal consistency across any separable metric spaces. Additionally, we propose another estimator for the conditional Fréchet mean, leveraging data-driven partitioning and quantization, and demonstrate its universal consistency when the output space is any Banach space. \\
\end{abstract}

\section{Introduction}

Let $(\X,\rho_{\X})$ and $(\Y,\rho_{\Y})$ be two separable metric spaces equipped with their Borel $\sigma$-algebra and let $(X, Y)$ be a pair of two random variables taking values in $\X$ and $\Y$, respectively. 
We denote by $\mu_X$ and $\mu_Y$ their respective marginal probability distributions.
In this work, we study the problem of estimating the Fréchet mean and the conditional Fréchet mean, when both variables $X$ and $Y$ take metric space values.

Given a loss function $\ell: \Y\times\Y\to \R_+$, if
\[
R(y) = \e{\ell(Y,y)}
\]
denotes the risk at $y$, then the {\it Fréchet mean} is defined by
\begin{equation}\label{eq:frech_mean}
    \bar{m}_{Fr} = \argmin_{y\in \Y} R(y) \; ,
\end{equation}
also called barycenter or centroid \citep{frechet1948elements}.
In this paper, we assume that the  Fréchet mean $\bar m_{Fr}$ is well defined, which implies that $\argmin_{y\in \Y}$ exists. However, in general, the uniqueness of the minimizer cannot be guaranteed, cf. \cite{hein2009robust}.
Instead, we assume that there is a (possibly unknown) measurable selection rule among the set of minimizers. \\

Next, we review some examples where our estimator of $\bar{m}_{Fr}$ can be applied.

\begin{example}
\label{sq}
Let $\Y$ be a separable Hilbert space with norm $\|\cdot \|$ and 
\[
\ell(y,y')=\|y-y'\|^2 \; .
\]
If $\Y=\R$, then this cost is called the squared loss.
Due to
\begin{align*}
\e{\ell(Y,y)}
&=
\e{\|Y-y\|^2}
=
\e{\|Y-\e{Y}\|^2}+\|\e{Y}-y\|^2 \; ,
\end{align*}
one gets that
\begin{align*}
\bar{m}_{Fr}
&=
\e{Y},
\end{align*}
that is, the Fréchet mean is equal to the expectation.
Concerning the definition of expectation in separable Banach space, we refer to \cite{beck1963}. 
The estimation problem of the expectation $\e{Y}$ is easy, because the sample mean is strongly consistent.
\end{example}

\begin{example}
\label{ab}
With respect to Example \ref{sq}, we modify the cost to:
\[
\ell(y,y')=\|y-y'\| \; .
\]
If $\Y=\R$, then this cost is called the absolute loss and $\bar{m}_{Fr}$ is the median, which may not be unique, cf. \cite{milasevic1987uniqueness}. 
In this setting, estimating the median is easy, since the sample median is consistent.
However, the setup of $\Y=\R^d$ with $d\geq 2$ is much more difficult since there is no known sample estimate of the corresponding Fréchet mean.
\end{example}

\begin{example}
\label{den}
Let  $\Y\subset L_1(\R^d)$ be the set of probability densities defined on $\R^d$, that is, $Y$ is a random density.
Given two densities $f$ and $g$, put
\[
\ell(f,g)
=\|f-g\|_1
=\int_{\R^d}|f(z)-g(z)|\lambda(dz) \; ,
\]
where $\lambda$ stands for the Lebesgue measure.
Then, $\bar{m}_{Fr}$ is the centroid of random densities.\\
In an equivalent setup, $\Y$ is the set of absolutely continuous probability measures on $\R^d$.
If $\nu_f$ and $\nu_g$ are probability measures with densities $f$ and $g$, respectively, then put
\[
\ell(\nu_f,\nu_g)
=TV(\nu_f,\nu_g)
\]
with total variation distance $TV$.
By Scheffé theorem,
\[
TV(\nu_f,\nu_g)
=\frac 12 \|f-g\|_1 \; .
\]
This implies that if $f^*$ is the centroid of densities and $\nu^*$ is the centroid of the corresponding probability measures, then
\[
\nu^*=\nu_{f^*} \; .
\]
\end{example}

\begin{example}
\label{gra}
In network analysis, $\Y$ can be a space of graphs or the space of graph Laplacians (see e.g. \citet{ginestet2012strong, ginestet2017hypothesis, kolaczyk2020averages, lunagomez2021modeling, severn2022manifold}).
\end{example}

\begin{example}
\label{fu}
\citet{bhattacharya2003large, bhattacharya2005large} study the case where $\Y$ is a Riemannian manifold.
\citet{ahidar2020convergence} provided convergence rates in compact geodesic metric spaces under general conditions.
Recently, \citet{lugosi2024uncertainty} demonstrated some practical examples of clinical data like MRI images or biosensor curves.
\end{example}

For the regression problem with metric space values, $X$ is called predictor, while $Y$ is the corresponding response.
We are looking for a function $m:\X\to \Y$, for which
$m(X)\approx Y$. This approximation is qualified by the loss function $\ell$
such that the risk of $m$ is defined by
\begin{align}
\label{risk}
R(m)=\e{\ell(Y,m(X))} \; .
\end{align}
In the standard regression setup, both $\X$ and $\Y$ are finite-dimensional Euclidean spaces and $\ell(y,y')=\|y-y'\|^2$, where $\|\cdot \|$ stands for the Euclidean norm. 
In this article, we consider a more general setup as $X$ and $Y$ can take values in more complex spaces such as spaces of curves or words or other infinite-dimensional spaces such as functional spaces. 
The best possible regression function is then assumed to be the one minimizing the risk in equation ~\eqref{risk}, and its solution is given by the regression function:
\begin{align}
\label{regr}
m_{Fr}(x)=\argmin_{y\in \Y} \e{\ell(Y,y)\mid X=x} \; ,
\end{align}
where $x\in \X$. The function $m_{Fr}$ is often referred to as the {\it conditional Fréchet mean} of $Y$ given $X$ or the Fréchet regression function.

\begin{example}
\label{sqX}
In the setup of Example \ref{sq} one gets that
\begin{align*}
\e{\ell(Y,y)\mid X=x}
&=
\e{\|Y-y\|^2\mid X=x}\\
&=
\e{\|Y-\e{Y\mid X=x}\|^2\mid X=x}+\|\e{Y\mid X=x}-y\|^2 \; .
\end{align*}
Therefore,
\begin{align*}
m_{Fr}(x)
&=
\e{Y\mid X=x} \; ,
\end{align*}
that is, the conditional Fréchet mean $m_{Fr}(x)$ is equal to the conditional expectation\\ $\e{Y\mid X=x}$.
If $\X$ is any separable metric space and $\Y=\R$, then \cite{hanneke2021universal} and  \cite{gyorfi2021universal} introduced universally consistent estimates of $m_{Fr}(\cdot)$.
\end{example}

\begin{example}
\label{abX} 
For the loss of Example 2, 
the conditional Fréchet mean is an extension of quantile regression.
For the particular case $\X=\R^d$ and $\Y=\R$, we refer to \cite{koenker2001quantile}.
\end{example}

\begin{example}
\label{denX} 
Recently, \citet{nguyen2026bayesian}
introduced an interesting setup of density-density regression scheme such that 
both $\X$ and $\Y$ are the sets of multivariate continuous distribution functions where the metrics are the Wasserstein metrics.
\end{example}

\begin{example}
\label{fuX}
If $\X=\R^d$ and $\Y$ is the space of Graph Laplacians, then \citet{zhou2022network} studied the estimate of the corresponding conditional Fréchet mean.
\end{example}

In the existing literature, the Fréchet mean or conditional mean estimators are obtained by minimizing a loss function on $\Y$. 
However, in practice, computing this minimum on a metric space is not always straightforward. In the Fréchet mean case, one possibility is to use what is called the restricted Frechet means where the minimum is not taken over $\Y$ but only over a sample (see e.g. \citet{sverdrup1981strong, ginestet2012strong}). In the conditional Fréchet mean case, \cite{petersen2019frechet} proposed a computational method for $\Y$ being a specific metric space, but no systematic computational method was given.
Although this issue highlights the complexity of practical implementation for general $\Y$ or general loss, recent advances have attempted to address it. 
In particular, \cite{cohen2022learning} developed a universally consistent Fréchet regression estimate together with an algorithm. 
Their approach is to first construct $\delta$-nets for $(X_1, \ldots, X_n)$ across multiple values of $\delta$ and then, to output the medoid point of each resulting cell through a particular minimization problem. 
However, a significant limitation of their method lies in the fact that when $\Y$ is not a countable metric space, they assume access to an oracle process capable of providing an $\varepsilon$-net of $\Y$, which is a rather strong assumption. 
Overall, in general, minimization is typically addressed on a case-by-case basis, depending on the specific structure of the metric space under consideration.
The main aim of our paper is therefore to come up with a universally consistent Fréchet regression estimate with small computational complexity. Our contributions can be summarized as follows.

\paragraph{Contributions.} 
\begin{enumerate}
    \item In Section \ref{sec:Frechet_mean}, a novel practical estimate of the Fréchet mean is introduced, differing slightly in its construction from the restricted Fréchet mean of \citet{sverdrup1981strong}. Under a boundedness condition on the loss function, we demonstrate that the estimator is strongly consistency, as soon as the Fréchet mean belongs to the support of the distribution of $Y$, Theorem \ref{thm:mean}.
    \item Concerning the conditional Fréchet mean, in Section \ref{sec:Frechet_reg}, we introduce an extension of the {\sc Proto-NN} estimate from \cite{gyorfi2021universal}, when both $\X$ and $\Y$ are separable metric spaces. We prove that, under mild assumptions, our estimate is strongly universally consistent, i.e., the risk evaluated at the estimate converges almost surely to the risk of the conditional Fréchet mean when it belongs to the support of the distribution of $Y$, Theorem \ref{thm:regr}. 
\end{enumerate}

\section{Empirical Fréchet mean}\label{sec:Frechet_mean}

Given some data $Y_1,\dots, Y_n$ from a distribution $\mu_Y$, a natural estimate of the Fréchet mean (Equation \eqref{eq:frech_mean}) is
\begin{align}
\label{mean}
\bar{m}_{Fr,n}=\argmin_{y\in \Y}  \sum_{i=1}^n{\ell(Y_i,y)} \; .
\end{align}
The asymptotic theory of the empirical Fréchet mean is now well established. For example, when the space $\Y$ is a separable finite quasi-metric space, its large-sample properties were analyzed by \citet{ziezold1977expected}. 
Beyond asymptotic theory, non-asymptotic results have also been established in various settings.  \citet{le2022fast} and \citet{schotz2019convergence} prove upper-bounds on the distance between the Fréchet mean and its empirical counterpart. 
 We also refer to \cite{schotz2022strong, aveni2024uniform, evans2024limit, hotz2024central, jaffe2024fr} and the references therein for recent results.
More recently, \citet{brunel2024concentration} gave finite sample concentration inequalities for the same quantity. 
Finally, significant progress has also been made in the study of robust Fréchet mean, with notable contributions including \cite{hsu2016loss, yun2023exponential, kim2025robust} and \cite{bartl2025robust}.
Note that these papers do not consider the difficult algorithmic question of how to solve the minimization over the full space $\Y$ in practice. \\

To our knowledge, the only estimate together with an algorithm for general space $\Y$ has been defined by \citet{sverdrup1981strong} as follows:
\begin{align}
\label{ST}
\widetilde{m}_{n}=\argmin_{y \in \{Y_j\}_{j=1}^n} \sum_{i=1}^n{\ell(Y_i,y)} \; .
\end{align}
This estimate is called the restricted Fréchet means and can be interpreted as a resubstitution estimate. Its strong consistency is proved under the condition that the Fréchet mean belongs to the support of $Y$, see condition \eqref{mSY} below. \\

\noindent
Our estimate is as follows:
\paragraph{Algorithm 1.} As above, the idea is to quantize the space $\Y$, making the minimization possible. To this aim, we assume to have access to a set of i.i.d. copies $Y'_1,\dots ,Y'_n$ of the label $Y$.
Let $\Y_n$ be the set of prototypes $Y'_1,\dots ,Y'_n$, our Fréchet mean estimate is defined by
\begin{align}
\label{emean}
\bar{m}_{n}=\argmin_{y\in \Y_n}  \sum_{i=1}^n{\ell(Y_i,y)} \; .
\end{align}

Our estimate can be interpreted as a splitting data estimate, where the $2n$ samples are splitted into the learning samples $Y_1,\dots ,Y_n$ and the testing data $Y'_1,\dots ,Y'_n$. Notice that, no matter the number of prototypes in $\Y_n$, $m_n$ takes values in the support $\suppY$ of $\mu_Y$ -- its existence being guaranteed because $\Y$ is separable (Theorem 2.1 p. 27 in \citet{parthasarathy2005probability}). 
On the theoretical side, both algorithms are universally consistent  under the same set of assumptions, but considering Algorithm 1 the proof is simpler. \\

\noindent
\textbf{Assumption 1}
\textit{
If $\suppY\varsubsetneq \Y$, this can raise consistency issues as $\bar m_{Fr}$ can take values in $\Y\backslash \suppY$ while $\bar m_{n}$ will never.  Therefore, in the following, we suppose that:
\begin{align}
\label{mSY}
\bar m_{Fr}\in \suppY\; .
\end{align} 
}

This assumption implies that
\begin{align*}
\min_{y\in \Y} R(y)
&=
\min_{y\in \suppY} R(y) \; .
\end{align*}
We can now state our consistency result.
\begin{theorem}
\label{thm:mean}
Assume \eqref{mSY}.
Furthermore, suppose that for all $y,y',y''\in \Y$,
\[
\sup_{y,y'\in \Y}\ell(y,y')\leq L< \infty \; ,
\]
\begin{align}
\label{tr}
|\ell(y,y')-\ell(y,y'')|
&\le c\cdot \ell(y',y'')
\end{align}
for a finite constant $c$,
and 
\begin{align}
\label{lr}
\ell(y,y')
&\le  \rho_{\Y}(y,y')^{\alpha}
\end{align}
with $\alpha>0$. 
Then
\begin{align}
\label{mapr}
\lim_n R(\bar m_n)=R(\bar m_{Fr}) \quad a.s.
\end{align}
\end{theorem}

\section{Fréchet regression estimate}\label{sec:Frechet_reg}

In order to estimate the Fréchet regression function $m_{Fr}$ (Equation \ref{regr}), we are given a set of $n$ labeled training data:
\begin{align}
\label{data}
\D_n=\{(X_1,Y_1), \dots ,(X_n,Y_n)\}
\end{align}
all drawn i.i.d. from the same distribution. In this section, we are interested in constructing an estimator $\text{Est}_n=\text{Est}_n(\D_n):\X\rightarrow \Y$ of $m_{Fr}$ that is \textit{universally strongly consistent}. In other words, we look for $\text{Est}_n$ such that $\PROB\left\{ \underset{n\rightarrow \infty}{\lim}R(\text{Est}_n)=R(m_{Fr})\right\}=1$, for every distribution of $(X,Y)$. \\

In the past, the Fréchet regression estimation problem has been studied mainly when either the predictor space $\X$ or the response space $\Y$ was non-Euclidean.
When $\X$ is a separable metric space and $\Y=\R$, the most popular estimator is the k-Nearest-Neighbor (k-NN) estimator.
Unfortunately, the k-NN estimate is not universally consistent, see
\cite{cerou2006nearest}, \cite{collins2020universal} and \cite{kumari2024universal}. 
The estimate  {\sc OptiNet} in \cite{hanneke2021universal} and the estimate  {\sc Proto-NN} in \cite{gyorfi2021universal} are universally consistent regression estimates, see the extension of {\sc Proto-NN} in the next section. 

In general, relatively little attention has been paid to the more general case where both $\X$ and $\Y$ are metric spaces.
For example, in \cite{hein2009robust} and \cite{steinke2010nonparametric}, the authors introduce a Nadaraya-Watson based estimator, but the consistency of this estimator is only established under the restrictive assumption that $\X$ and $\Y$ are both Riemannian manifolds. 
In \cite{petersen2019frechet}, the authors extended the local linear regression in \cite{fan1996local} to take into account that $\Y$ is a metric space, but their framework assumes that $\X$ is $\R^d$. Furthermore, their proof of pointwise consistency relies on specific distribution assumptions, which excludes universal consistency. 
Although \citet{chen2022uniform} subsequently improved these results by providing uniform convergence rates, their conclusions still depend on distributional assumptions. More recently, \citet{capitaine2024frechet} introduced Fréchet trees and Fréchet random forests, and \citet{qiu2024random} proposed a random forest-weighted local Fréchet regression method for cases where both $\X$ and $\Y$ are metric spaces. However, their theoretical results are again based on assumptions on the distributions. \\

\paragraph{Construction of our regression estimator:}
We now present our universally consistent estimator of the Fréchet regression function $m_{Fr}$, which is a natural extension of the {\sc Proto-NN} estimate introduced in \citet{gyorfi2021universal}. The main difference between ours and the one in \citet{gyorfi2021universal} is the fact that in the former paper the output space $\Y$ is a set of real values, while in the present we consider $\Y$ to be any separable Banach space. \\

The estimator construction works as follows. For an integer $k\geq 1$, we assume that in addition to the labeled sample $\D_n$, we also have access to an independent set of unlabeled samples, called prototypes and denoted by
$\X_k =\{X'_1,\dots ,X'_k\}$, where the $X'_i$'s are independent copies of $X$.
Let the data-driven partition $\P_{k}$ of $\X$ be such that $\P_{k}$ is a Voronoi partition with the nucleus set $\X_k$, i.e., 
\begin{align}
\label{Pk}
\P_{k}=\{A_{k,1},A_{k,2},\dots ,A_{k,k}\}
\end{align}
such that $A_{k,j}$ is the Voronoi cell around the nucleus $X'_j$,
\begin{align*}
A_{k,j} = \left\{x\in\X : j = \argmin_{1\leq i \leq k} \rho_{\X}(x,X'_i)\right\},
\end{align*}
where tie breaking is done by indices, i.e., if $X'_i$ and $X'_j$ are equidistant from $x$, then $X'_i$ is declared ``closer'' if $i < j$.
In this paper, we assume that ties occur with probability zero.

Given that, we are able to define an estimator for the conditional expectation\\
$\e{\ell(Y,y)\mid X=x}$ that we wish to minimize over $\Y$ (see Equation \ref{regr}). While a naive estimator could be the empirical one
\begin{equation*}
    \frac{\sum_{i=1}^n\ell(Y_i,y)\IND_{\{X_i = x \}}}{\sum_{i=1}^n\IND_{\{X_i = x \}}}\;,
\end{equation*}
it must be noted that for many data-distributions (e.g. any continuous random variable $X$), the event $\{X=x\}$ can be of null probability measure, leading to a poor, or non-computable, estimate. Instead, we propose to approximate the conditional expectation $\e{\ell(Y,y)\mid X=x}$ by
$\e{\ell(Y,y)\mid X\in A_{k,j}}$ for all $x\in A_{k,j}$. In the end, for $x\in A_{k,j}$, we take:
\begin{equation*}
    \frac{\sum_{i=1}^n\ell(Y_i,y)\IND_{\{X_i\in A_{k,j} \}}}{\sum_{i=1}^n\IND_{\{X_i\in A_{k,j} \}}}\;,
\end{equation*}
to be our final estimate of the conditional expectation in Equation \eqref{regr}. A natural way to estimate $m_{Fr}(x)$ for $x\in A_{k,j}$ would be to find $y\in\Y$ that minimizes the above empirical estimate.
In general, the minimization over the whole space $y\in\Y$ is impossible,
therefore we present an algorithm (referred to as Algorithm 2), which is a minimization over an appropriate quantization of the space, analogous to the one of Algorithm 1. 

\paragraph{Algorithm 2.} As in the previous section, let $\Y_n$ be the set of prototypes $Y'_1,\dots ,Y'_n$. Then our Fréchet regression estimate is defined by
\begin{align*}
m_{n}(x)=\argmin_{y\in \Y_n} \sum_{i=1}^n\ell(Y_i,y)\IND_{\{X_i\in A_{k,j} \}}, \quad \mbox{ if } 
x\in A_{k,j}.
\end{align*}

Hence, the estimate $m_{n}(\cdot)$ is a piece-wise constant function, defined on Voronoi cells $A_{k,1},\ldots, A_{k,k}$, and with values taken in $\Y_n$. Notice that, no matter the number of prototypes in $\Y_n$, $m_n$ takes values in the support $\suppY$ of $\mu_Y$ -- its existence being guaranteed because $\Y$ is separable (Theorem 2.1 p. 27 in \citet{parthasarathy2005probability}). \\

\noindent
\textbf{Assumption 2}
\textit{
As for the Fréchet mean, if $\suppY\varsubsetneq \Y$, this can raise consistency issues as $m_{Fr}(x)$ can take values in $\Y\backslash \suppY$ while $m_{n}(x)$ will never.  Therefore, in the following we assume that for all $x\in \X$, 
\begin{align}
\label{SY}
m_{Fr}(x)\in \suppY\; .
\end{align} 
}
It is a mild condition.
For example, if $\Y$ is a vector space and $\suppY$ is convex, then this condition is satisfied.
The condition \eqref{SY} implies that 
\begin{align*}
\min_{y\in \Y} \e{\ell(Y,y)\mid X=x}
&=
\min_{y\in \suppY} \e{\ell(Y,y)\mid X=x} \; .
\end{align*}

We can now state our main consistency result.
\begin{theorem}
\label{thm:regr}
Assume \eqref{SY}, that $\Y$ is a separable Banach space, and that  for $\X_k$, ties occur with probability zero.
Suppose that the conditions on the cost $\ell$ of Theorem \ref{thm:mean} are satisfied.
If $k=k_n$ is such that $k_n\to \infty$ and $k_n\ln n/n\to 0$,
then
\begin{align}
\label{apr}
\lim_n R(m_n)=R(m_{Fr}) \quad \; a.s.
\end{align}
\end{theorem}

The limit relation~\eqref{apr} in Theorem \ref{thm:regr} provides the strong consistency guarantee we were looking for. 
It is universal because no assumption is made regarding the distributions of the data and the Fréchet regression function $m_{Fr}(\cdot)$. 
Interestingly, consistency holds even in the difficult case, when the risk $R(m)$ has several minima.
The only assumptions are with respect to the loss function $\ell$. 
The boundedness assumption with $L$ is rather standard, it is, for instance, referred to as the bounded diameter assumption in \cite{cohen2022learning}. 
As a matter of fact, when $\ell$ is a metric the other two assumptions in equations \eqref{tr} and \eqref{lr} are not necessary. Indeed, \eqref{tr} would be directly satisfied with $c=1$ by the reverse triangle inequality, same for \eqref{lr} since it is here only to ensure that $\ell$ is uniformly bounded by a continuous function.

\section{Proofs} \label{sec:proofs}

\subsection{Proof of Theorem \ref{thm:mean}}

Note that $\Y_{n}$ is a finite set of size $n$ and $\bar m_n$ takes values in $\Y_{n}$, i.e., if we denote the empirical risk by 
\begin{align*}
R_n(y) := \frac{1}{n} \sum^n_{i=1} \ell(y, Y_i) \, ,
\end{align*}
then,
\begin{align}
\label{mmn}
\bar m_n=\argmin_{y \in \Y_{n} }R_n(y) \; .
\end{align}

We now decompose the excess risk in two terms, namely, the estimation error and the approximation error:
\begin{align*}
    R(\bar m_n) - R (\bar m_{Fr})
    &= R(\bar m_n) - \min_{y \in \Y_{n}} R(y)
    + \min_{y \in \Y_{n}} R(y) - R(\bar m_{Fr}) \; .
\end{align*}

\subsection*{Controlling the estimation error $R(\bar m_n) - \min_{y \in \Y_{n}} R(y)$}

The {\it estimation error} is bounded as follow:
\begin{align*}
     R(\bar m_n) - \min_{y \in \Y_{n}} R(y) 
    &= \max_{y \in \Y_{n}} \paren{R(\bar m_n) - R(y)} \\
    &= \max_{y \in \Y_{n}} \paren{R(\bar m_n) - R_n(\bar m_n) + R_n(\bar m_n) - R_n(y) + R_n(y) - R(y)} \\
    &\leq \max_{y \in \Y_{n}} \paren{R (\bar m_n) - R_n(\bar m_n) + R_n(y) - R(y)} \\
    &\leq 2 \max_{y \in \Y_{n}} \abs{R_n(y) - R(y)} \; ,
\end{align*}
where we have used that for any $y \in \Y_{n}$, by the definition (\ref{mmn}),
$R_n(\bar m_n) - R_n(y) \leq 0$. 
Finally, the union bound together with the Hoeffding inequality and the assumption that $\ell(y, y') \leq L$ imply that, for all $\varepsilon > 0$:
\begin{align*}
    \probp{\max_{y \in \Y_{n}} \abs{R_n(y) - R(y)} > \varepsilon \mid \Y_n} 
    &\leq n \max_{y \in \Y_{n}}\probp{\abs{R_n(y) - R(y)} > \varepsilon}\\
    &\leq 2 n e^{-2 n \varepsilon^2 / L^2}\\
    &= 2 e^{-2 n \varepsilon^2 / L^2+ \ln n}\; .
\end{align*}
Therefore,  the Borel-Cantelli lemma implies that
\begin{align*}
    R(\bar m_n) - \min_{y \in \Y_{n}} R(y) \rightarrow 0 \quad a.s.
\end{align*}
and the consistency of the estimation error is proved. \\

\subsection*{Controlling the approximation error $\min_{y \in \Y_{n}} R(y) - R(\bar m_{Fr})$}

By condition (Assumption 1), the approximation error $I_n$ has the form
\begin{align*}
    I_{n}
    &= R(\bar m_{n})-R(\bar m_{Fr}) 
    = \min_{y \in \Y_{n}} R(y)-\min_{y \in \Y} R(y) 
    = \min_{y \in \Y_{n}} R(y)-\min_{y \in \suppY} R(y).
\end{align*}
The task is to show that 
\begin{align*}
I_{n} & = \min_{y \in \Y_{n}} R(y)-\min_{y \in \suppY} R(y)
\rightarrow 0 \quad \text{a.s.}
\end{align*}
For any  $0<s$, one gets
\begin{align*}
	\PROB\{I_{n}>s \} 
    &=\PROB\left\{ \min_{y \in \Y_{n}} R(y)-\min_{y \in \suppY} R(y)  > s \right\}\\
    &=\PROB\left\{ \min_{i=1,\dots , n} R(Y'_i)-\min_{y \in \suppY} R(y)  > s \right\}\\    
    &=\PROB\left\{ R(Y'_1)-\min_{y \in \suppY} R(y)  > s \right\}^n .
\end{align*}
Therefore,
\begin{align*}
\sum_{n=1}^\infty \PROB\{I_{n} > s\}
&<\frac{1}{\PROB\left\{ R(Y'_1)-\min_{y \in \suppY} R(y)  \le  s \right\}\}}
<\infty ,
\end{align*}
if
\begin{align*}
\PROB\left\{ R(Y'_1)-\min_{y \in \suppY} R(y)  \le  s \right\}
&> 0.
\end{align*}
Let $\Tilde{y} = \arg\min_{y \in \suppY} R(y)$,
then
\begin{align*}
    \PROB\left\{ R(Y'_1)-\min_{y \in \suppY} R(y)  \le  s \right\}
    &=\PROB\left\{\EXP[\ell(Y, Y'_1) \mid Y'_1] - \EXP[\ell(Y, \Tilde{y})]  \le s \right\} \\
    &\ge \PROB\left\{c\cdot \ell(Y'_1,\Tilde{y})  \le s \right\} \\
    &\geq \PROB\left\{\rho_{\Y}(Y'_1, \Tilde{y})  \le (s/c)^{1/\alpha} \right\}\\
    &> 0
\end{align*}
where we have used \eqref{tr} and \eqref{lr} for the first two inequalities and the last inequality holds because of the property of the support $\suppY$. This ends the proof. 

\hfill $\square$

\subsection{Proof of Theorem \ref{thm:regr}}

Let us first introduce the set of functions that take their values in $\cY_n$ and that are constant over the Voronois cells $\set{A_{k, j}}^k_{j=1}$:
\begin{align*}
    \cF_{n, k} := \cF(\cX_k, \cY_n) 
    := \set{f :  f(x) = \sum^k_{j=1} y_j \IND\set{x \in A_{k, j}}, y_j \in \cY_n } \; .
\end{align*}

Note that we have $k$ distinct Voronoi cells because we assume that there is no tie. Hence,
$\cF_{n, k}$ is a finite set of size $n^k$ and $m_n$ belongs to $\cF_{n, k}$, i.e., if we denote the empirical risk by 
\begin{align*}
R_n(f) := \frac{1}{n} \sum^n_{i=1} \ell(f(X_i), Y_i) \, ,
\end{align*}
then,
\begin{align}
\label{mn}
m_n=\argmin_{f \in \cF_{n, k} }R_n(f) \; .
\end{align}

We now decompose the excess risk in two terms, namely, the estimation error and the approximation error:
\begin{align*}
    R(m_n) - R(m_{Fr})
    &= R(m_n) - \min_{f \in \cF_{n, k}} R(f)
    + \min_{f \in \cF_{n, k}} R(f) - R(m_{Fr}) \; .
\end{align*}
The proof is divided into two parts. One part to show that the estimation error\\
$R(m_n) - \min_{f \in \cF_{n, k}} R(f)$ tends to $0$ a.s., and the other part is showing that the approximation error $\min_{f \in \cF_{n, k}} R(f) - R(m_{Fr})$ also tends to $0$ a.s. 
The control of the approximation error is also split in two terms, because, as explained later, two types of approximation errors appear: one due to the discretizations of the $\Y$ space with the prototypes in $\Y_n$, and one due to the approximation of the $\X$ space with the Voronoi partition.

\subsection*{Controlling the estimation error $R(m_n) - \min_{f \in \cF_{n, k}} R(f)$}
The {\it estimation error} is bounded as follow:
\begin{align*}
    & R(m_n) - \min_{f \in \cF_{n, k}} R(f) \\
    &= \max_{f \in \cF_{n, k}} \paren{R(m_n) - R(f)} \\
    &= \max_{f \in \cF_{n, k}} \paren{R(m_n) - R_n(m_n) + R_n(m_n) - R_n(f) + R_n(f) - R(f)} \\
    &\leq \max_{f \in \cF_{n, k}} \paren{R(m_n) - R_n(m_n) + R_n(f) - R(f)} \\
    &\leq 2 \max_{f \in \cF_{n, k}} \abs{R_n(f) - R(f)} \; ,
\end{align*}
where we have used that for any $f \in \cF_{n, k}$, by the definition (\ref{mn}),
$R_n(m_n) - R_n(f) \leq 0$. 
Finally, the union bound, the Hoeffding inequality and the assumption that $\ell(y, y') \leq L$ imply that, for all $\varepsilon > 0$:
\begin{align*}
    \probp{\max_{f \in \cF_{n, k}} \abs{R_n(f) - R(f)} > \varepsilon \mid \X_k, \Y_n} 
    &\leq \lvert \cF_{n, k} \rvert  \cdot \max_{f \in \cF_{n, k}}\probp{\abs{R_n(f) - R(f)} > \varepsilon}\\
    &\leq 2 n^k e^{-2 n \varepsilon^2 / L^2}\\
    &= 2 e^{-2 n \varepsilon^2 / L^2+k \ln n}\; .
\end{align*}
Therefore, the condition $k_n\ln n/n\to 0 $ together with the Borel-Cantelli lemma implies that
\begin{align*}
    R(m_n) - \min_{f \in \cF_{n, k}} R(f) \rightarrow 0 \quad a.s.
\end{align*}
and the consistency of the estimation error is proved. \\

\subsection*{Controlling the approximation error $\min_{f \in \cF_{n, k}} R(f) - R(m_{Fr})$.}

First, remark that, given a Voronois cell $A_{k_n, j}$, when a new point $X_{k_n + 1}'$ is observed, the cell $A_{k_n, j}$ either remains unchanged or is partitioned into several new Voronoi cells. This means that $\cF_{n, k_n} \subseteq \cF_{n+1, k_{n}}$, $\cF_{n, k_n} \subseteq \cF_{n, k_{n+1}}$, and $\cF_{n, k_n} \subseteq \cF_{n+1, k_{n+1}}$. 
Therefore, for any fixed $k$, 
\[
\min_{f \in \cF_{n, k}} R(f)\downarrow
\]
as $n\to \infty$, and for any fixed $n$, 
\[
\min_{f \in \cF_{n, k}} R(f)\downarrow
\]
as $k\to \infty$.
It implies that for any $k_n\uparrow$, 
\begin{align}
\label{Dow}
\min_{f \in \cF_{n, k_n}} R(f)\downarrow
\end{align}
as $n\to \infty$. \\

We have to show that
\begin{align}
\label{AS}
\lim_n \left(\min_{f \in \cF_{n, k_n}} R(f) -R(m_{Fr})\right)=0 \qquad a.s. 
\end{align}
Because of  (\ref{Dow}), we know that $\min_{f \in \cF_{n, k_n}} R(f) -R(m_{Fr})$ almost surely converges to a non-negative random variable. Since a non-negative random variable with expectation equal to $0$ is almost surely $0$, to show~\eqref{AS} it is sufficient to show that the expectation of this limit is equal to $0$.  By the dominated convergence theorem, we have:
\begin{align*}
\EXP\left\{\lim_n \left(\min_{f \in \cF_{n, k_n}} R(f) -R(m_{Fr})\right)\right\}
&=
\lim_n \left(\EXP\left\{\min_{f \in \cF_{n, k_n}} R(f)\right\} -R(m_{Fr})\right) \; .
\end{align*}
Therefore, the task (\ref{AS}) is reduced to show that 
\begin{align}
\label{DS}
\lim_n \left(\EXP\left\{\min_{f \in \cF_{n, k_n}} R(f)\right\} -R(m_{Fr})\right)=0 \; .
\end{align}

One gets that
\begin{align*}
\lim_n \EXP\left\{\min_{f \in \cF_{n, k_n}} R(f)\right\}
&\ge
\inf_{k\in\mathbb{N}} \inf_{n \in\mathbb{N}} \EXP\left\{\min_{f \in \cF_{n, k}} R(f)\right\}.
\end{align*}
For any fixed $k$, the monotonicity property implies that
\begin{align*}
\lim_n \EXP\left\{\min_{f \in \cF_{n, k_n}} R(f)\right\}
&\le
\lim_n \EXP\left\{\min_{f \in \cF_{n, k}} R(f)\right\}
=
\inf_n \EXP\left\{\min_{f \in \cF_{n, k}} R(f)\right\}.
\end{align*}
Therefore,
\begin{align*}
\lim_n \EXP\left\{\min_{f \in \cF_{n, k_n}} R(f)\right\}
&\le
\inf_k\inf_n \EXP\left\{\min_{f \in \cF_{n, k}} R(f)\right\}.
\end{align*}
Thus,
\begin{align}
\label{ii}
\lim_n \EXP\left\{\min_{f \in \cF_{n, k_n}} R(f)\right\}
&=
\inf_k \inf_n \EXP\left\{\min_{f \in \cF_{n, k}} R(f)\right\}.
\end{align}

This equality together with (\ref{DS}) yields that the only task is to show that
\begin{align}
\label{iii}
\lim_k \lim_n \EXP\left\{\min_{f \in \cF_{n, k}} R(f)\right\}
-R(m_{Fr})
&=
0.
\end{align}

Now, let us introduce 
\begin{align*}
    \cF_{\infty, k} := \set{f :  f(x) = \sum^k_{j=1} y_j \IND\set{x \in A_{k, j}}, y_j\in \suppY } \; .
\end{align*}
One gets that
\begin{align*}
 \EXP\left\{\min_{f \in \cF_{n, k}} R(f)\right\}-R(m_{Fr})
&=
I_{n,k}+J_{k}
\end{align*}
where
\begin{align*}
I_{n,k}
&=
\EXP\left\{\min_{f \in \cF_{n, k}} R(f)\right\} - \EXP\left\{\min_{f \in \cF_{\infty, k}} R(f)\right\}
\end{align*}
and
\begin{align*}
J_{k}
&=
\EXP\left\{\min_{f \in \cF_{\infty, k}} R(f)\right\}
-R(m_{Fr}) \; .
\end{align*}

\paragraph{Control of $I_{n, k}$:\\}
Recall that $\Y_n := \set{Y'_1, \ldots, Y'_n}$ and denote 
\begin{align*}
y^*_{k, j}=\displaystyle\argmin_{y \in \suppY} \e{\ell(Y, y) \IND\{X \in A_{k, j}\}}.
\end{align*}
Put
\begin{align*}
    G(y,x) & = \EXP\{\ell(Y,y)\mid X=x\} \; .
\end{align*}
For a function $f(x) = \sum^k_{j=1} y_j \IND\set{x \in A_{k, j}}$,
\begin{align*}
    R(f) 
    & = \int_{\X}\EXP\{\ell(Y,f(x))\mid X=x\}\mu(dx)\\
    &= \int_{\X}G(f(x),x)\mu(dx)\\
    &= \sum_{j=1}^k \int_{A_{k,j}}G(y_j,x)\mu(dx) \; .
\end{align*}
Then
\begin{align*}
    &\min_{f \in \cF_{n, k}} R(f) - \min_{f \in \cF_{\infty, k}} R(f) \\
    &=\min_{f \in \cF_{n, k}} \e{\ell(Y, f(X))} - \min_{f \in \cF_{\infty, k}} \e{\ell(Y, f(X))} \\
    &= \sum_{j=1}^k \min_{y \in \Y_n} \int_{A_{k,j}}G(y,x)\mu(dx)
    - \sum_{j=1}^k \min_{y \in \suppY} \int_{A_{k,j}}G(y,x)\mu(dx) \\
    &= \sum_{j=1}^k \paren{\min_{y \in \Y_n} \int_{A_{k,j}}G(y,x)\mu(dx) 
    - \int_{A_{k,j}}G(y^*_{k, j},x)\mu(dx)} \; ,
\end{align*}
which together with the conditions of the theorem implies that
\begin{align*}
    I_{n,k}
    &\leq c\sum_{j=1}^k \e{\min_{y \in \Y_n} \ell(y, y^*_{k, j}) \mu(A_{k,j})} \\
    &= c k \e{\min_{y \in \Y_n} \ell(y, y^*_{k, 1}) \mu(A_{k,1})} \\
    &\leq  c k \e{\min\{L,\min_{y\in \Y_n}\rho_\Y(y, y^*_{k, 1})^\alpha\}\mu(A_{k,1})} \\
    &= c k \e{\min\left\{L,\rho_\Y(Y'_{(1, n)}(y^*_{k, 1}), y^*_{k, 1})^\alpha \right\}\mu(A_{k,1})} \; ,
\end{align*}
where $Y'_{(1, n)}(y^*_{k,1})$ is the nearest neighbor of $y^*_{k, 1}$ in $\Y_n$. Because $\Y_n$ is a separable metric space and $y^*_{k, 1} \in \suppY$, using the Cover-Hart lemma \citep{cover1967nearest}, we know that $\rho_\Y(Y'_{(1, n)}(y^*_{k, 1}), y^*_{k, 1}) \rightarrow 0$ a.s. Using the dominated convergence theorem, this means that for any fixed $k$
\begin{align*}
    \lim_n I_{n,k}
    &=0 \; ,
\end{align*}
and therefore $\lim_k\lim_n I_{n,k}=0$ a.s.

\paragraph{Control of $J_{k}$:\\}
By condition, $\Y$ is a separable Banach space, therefore, combining the standard Lusin theorem \citep[Theorem 7.5.2]{dudley2018real} with the Dugundji extension theorem \citep{dugundji1951extension}, we can refer to a general version of the Lusin theorem as follows:
For any $\delta >0$, there exist a continuous function $\widehat{m}:\X\rightarrow\Y$ and a closed set $K \subset \cX$ such that $\mu(K^c) < \delta$, with $K^c := \cX \backslash K$,
and such that the functions $\widehat{m}$ and $m_{Fr}$ coincide on $K$. In fact, Lusin theorem gives the continuous function $\widehat{m}$ from $K$ to $\Y$. The Dugundji theorem allows to extend $\widehat{m}$ as a continuous function going from $\X$ to $\Y$ if $\Y$ is a locally convex topological vector space, which is the case of Banach spaces.\\

For the notation
\begin{align*}
f_k^*(x)
&= 
\widehat{m}(X'_{(1,k)}(x))
=
\widehat{m}(X'_j), \quad \mbox{if}\quad x\in A_{k,j},
\end{align*}
one gets
\begin{align*}
J_{k}
&=
\EXP\left\{\min_{f \in \cF_{\infty, k}} R(f)\right\}-R(m_{Fr})\\
&\le 
\EXP\left\{R(f_k^*)\right\}-R(m_{Fr})\\
&=
\EXP\left\{ \int  G(\widehat{m}(X'_{(1,k)}(x)),x)\mu(dx)\right\}-\int  G(m_{Fr}(x),x)\mu(dx)\\
&\le 
\left|\EXP\left\{ \int  G(\widehat{m}(X'_{(1,k)}(x)),x)\mu(dx)\right\}
- \int  G(\widehat{m}(x),x)\mu(dx)\right|\\
& \quad +
\left| \int  G(\widehat{m}(x),x)\mu(dx)-\int  G(m_{Fr}(x),x)\mu(dx) \right|.
\end{align*}

The second term is easy to control: 
\begin{align*}
\left| \int  G(\widehat{m}(x),x)\mu(dx)-\int  G(m_{Fr}(x),x)\mu(dx) \right|
&\le 
c\int  \ell (\widehat{m}(x) ,m_{Fr}(x) )\mu(dx)\\
&=
c\int_{K^c}  \ell(\widehat{m}(x) ,m_{Fr}(x) )\mu(dx)\\
&\le 
cL\mu(K^c) \\
&\le 
cL\delta .
\end{align*}
For the first term, the Cover-Hart lemma and the continuity of $\widehat{m}$ together with the dominated convergence theorem imply that
\begin{align*}
&
\left| \EXP\left\{ \int  G(\widehat{m}(X'_{(1,k)}(x)),x)\mu(dx)\right\}
- \int  G(\widehat{m}(x),x)\mu(dx)\right|\\
&\le 
c\int  \EXP\left\{\min\left\{L,\rho_\Y(\widehat{m}(X'_{(1,k)}(x)) , \widehat{m}(x))^\alpha\right\}\right\}\mu(dx) \\
&\to  
0 \; ,
\end{align*}
as $k\to \infty$.
Thus,
\begin{align*}
\limsup_k J_{k}
&\le
cL\delta \;.
\end{align*}
This concludes the proof.
\hfill $\square$

\section*{Acknowledgment}
P.H. gratefully acknowledges the Emergence project MARS of Sorbonne Université.

\vskip 0.2in
\bibliography{biblio}

\end{document}